\documentclass[pre,twocolumn,superscriptaddress]{revtex4}  

\usepackage[paperwidth=210mm,paperheight=297mm,centering,hmargin=2cm,vmargin=2.5cm]{geometry}

\usepackage[pdftex]{graphicx}
\usepackage{amsmath,amsfonts,amssymb}
\usepackage{morefloats}
\usepackage{xcolor} 
\usepackage{dsfont}

\usepackage{newfloat,algcompatible}
\usepackage{etoolbox}
\DeclareFloatingEnvironment[
    fileext=loa,
    listname=List of Algorithms,
    name= Algorithm,
    placement=tbhp,
]{algorithm}

\makeatletter
\renewcommand{\fnum@algorithm}{\small\textbf{\algorithmname~\thealgorithm}}
\makeatother

\definecolor{bluemoi}{rgb}{0.25,0.50 ,0.75} 

\usepackage{hyperref}
\hypersetup{
  bookmarksopen=false,
  pdftitle=" ",
  pdfauthor=" ", 
  pdfsubject=" ", 
  pdftoolbar=true, 
  pdfmenubar=true, 
  pdfhighlight=/O, 
  colorlinks=true, 
  pdfpagemode=UseNone, 
  pdfpagelayout=SinglePage, 
  pdffitwindow=true, 
  linkcolor=bluemoi, 
  citecolor=bluemoi, 
  urlcolor=bluemoi 
}

\makeatletter
\renewcommand{\figurename}{Figure}
\renewcommand{\fnum@figure}{\small\textbf{\figurename~\thefigure}}
\renewcommand{\thefigure}{\arabic{figure}}
\makeatother

\makeatletter
\renewcommand{\tablename}{Table}
\renewcommand{\fnum@table}{\small\textbf{\tablename~\thetable}}
\renewcommand{\thetable}{\arabic{table}}
\makeatother

\begin{document}

\title{Adaptive approximate Bayesian computation for complex models}

\author{Maxime Lenormand}\affiliation{IRSTEA, LISC, 24 avenue des Landais, 63172 AUBIERE, France}
\author{Franck jabot}\affiliation{IRSTEA, LISC, 24 avenue des Landais, 63172 AUBIERE, France}
\author{Guillaume Deffuant}\affiliation{IRSTEA, LISC, 24 avenue des Landais, 63172 AUBIERE, France}

\begin{abstract} 
We propose a new approximate Bayesian computation (ABC) algorithm that aims at minimizing the number of model runs for reaching a given quality of the posterior approximation. This algorithm automatically determines its sequence of tolerance levels and makes use of an easily interpretable stopping criterion. Moreover, it avoids the problem of particle duplication found when using a MCMC kernel. When applied to a toy example and to a complex social model, our algorithm is 2 to 8 times faster than the three main sequential ABC algorithms currently available.
\end{abstract}

\maketitle

Approximate Bayesian computation (ABC) techniques appear particularly relevant for calibrating stochastic models because they are easy to implement and applicable to any model. They generate a sample of model parameter values $(\theta_i)_{i = 1,..,N}$ (often also called particles) from the prior distribution $\pi(\theta)$ and select the $\theta_i$ values leading to model outputs $x \sim f(x|\theta_i)$ satisfying a proximity criterion with the target data $y$ ($\rho(x,y) \leq \epsilon$, $\rho$ expressing a distance, $\epsilon$ being a tolerance level). The selected sample of parameter values approximates the posterior distribution of parameters, leading to model outputs with the expected quality of approximation. However, in practise, running these techniques is very demanding computationally because sampling the whole space of parameters requires a number of simulations which grows exponentially with the number of parameters to identify. This tends to limit the application of these techniques to easily computable models \cite{Beaumont2010}. In this paper, our goal is  minimizing the number of model runs for reaching a given quality of posterior approximation, and thus to make the approach applicable to a larger set of models.

ABC is the subject of intense scientific researches and several improved versions of the original scheme are available, such as using local regressions to improve parameter inference \cite{Beaumont2002,Blum2010}, automatically selecting informative summary statistics \cite{Joyce2008,Fearnhead2011}, coupling to Markov chain Monte Carlo \cite{Marjoram2003,Wegmann2009} or improving sequentially the posterior distributions with sequential Monte Carlo methods \cite{Sisson2007,Toni2009,Beaumont2009}. This last class of methods approximates progressively the posterior, using sequential samples $S^{(t)}= (\theta_i^{(t)})_{i = 1,..,N}$ derived from sample $S^{(t-1)}$, and using a decreasing set of tolerance levels $\{\epsilon_1,...,\epsilon_T\}$. This strategy focuses the sampling effort in parts of the parameter space of high likelihood, avoiding to spend much computing time in systematically sampling the whole parameter space. 

The first sequential method applied to ABC was proposed by \cite{Sisson2007} with the ABC-PRC (Partial Rejection Control). This method is based on a theoretical work of \cite{DelMoral2006} to ABC. However, in \cite{Beaumont2009} the authors have shown that this method leads to a bias in the approximation of the posterior. In \cite{Beaumont2009,Toni2009} the authors proposed a new algorithm, called Population Monte Carlo ABC in \cite{Beaumont2009} and hereafter called PMC. This algorithm, corrects the bias by assigning to each particle a weight corresponding to the inverse of its importance in the sample. It is particularly interesting in our perspective because it provides with a rigorous framework to the sequential sample idea, which seems a good way for minimizing the number of runs. In this approach, the problem is then defining the sequence of tolerance levels $\{\epsilon_1,...,\epsilon_T\}$. In \cite{Drovandi2011} and \cite{DelMoral2012} the authors solve partly this problem by deriving the tolerance level at a given step from the previously selected sample. However, a difficulty remains: when to stop? If the final tolerance level $\epsilon_T$ is too large, the final posterior will be of bad quality. Inversely, a too small $\epsilon_T$ leads to a posterior that could have been obtained with less model runs.

In this paper, we propose a modification of the population Monte Carlo ABC algorithm proposed in \cite{Beaumont2009} that we call adaptive population Monte Carlo ABC (hereafter called APMC). This new algorithm determines by itself the sequence of tolerance levels as in \cite{Drovandi2011} and \cite{DelMoral2012}, and it also provides a stopping criterion. Furthermore, our approach avoids the problem of duplications of particles due to the MCMC kernel used in \cite{Drovandi2011} and \cite{DelMoral2012}. We prove that the computation of the weights associated to the particles in this algorithm lead to the intended posterior distribution and we also prove that the algorithm stops whatever the chosen value of the stopping parameter. We show that our algorithm, applied to a toy example and to an individual-based social model, requires significantly less simulations to reach a given quality level of the posterior distribution than the population Monte Carlo ABC algorithm of \cite{Beaumont2009} (hereafter called PMC), the replenishment SMC ABC algorithm of \cite{Drovandi2011} (hereafter called RSMC) and the adaptive SMC ABC algorithm of \cite{DelMoral2012} (hereafter called SMC). Our new algorithm has been implemented in the R package 'EasyABC' \cite{Jabot2013}. 

\section*{Sequential Monte-Carlo methods in approximate Bayesian computation}

In this section we present the three main sequential ABC algorithms currently available and their limitations. We present the Population Monte-Carlo ABC proposed in \cite{Beaumont2009} (hereafter called PMC), the Replenishment Sequential Monte-Carlo ABC proposed in \cite{Drovandi2011} and the Sequential Monte-Carlo ABC proposed in \cite{DelMoral2012}. These algorithms are detailed in Appendix A.

\subsection*{The PMC algorithm}

This method consists in generating a sample $S^{(t)}= (\theta_i^{(t)})_{i = 1,..,N}$ at each iteration of the algorithm, $1 \leq t \leq T$. Each particle of the sample $S^{(t)}$ satisfying the predefined tolerance level $\epsilon_t$ where $\epsilon_1 \geq \epsilon_t \geq \epsilon_T$. We say that a parameter value $\theta_i^{(t)}$, satisfies the tolerance level $\epsilon_t$, if when running the model we get $x \sim f(x|\theta_i^{(t)})$, such that its distance $\rho_i^{(t)}=\rho(x,y)$  to the target data $y$, is below $\epsilon_t$. At step $t$ the sample $S^{(t)}$ is derived from sample $S^{(t-1)}$ using a particle filter methodology. The first sample $S^1$ is generated using a regular ABC step. At step $t$ a new particle $\theta_i^{(t)}$ is generated using a Markov transition kernel $K_t$, $\theta_i^{(t)}\sim K_t(\theta|\theta*)$, until $\theta_i^{(t)}$ satisfies $\epsilon_t$ where $\theta^*$ is randomly draw from $S^{(t-1)}$ with probability $(w_i^{(t-1)})_{i = 1,..,N}$. The weight $w_i^{(t-1)}$ is proportional to the inverse of its importance in the sample $S^{(t-1)}$ (Eq. \ref{weightPMC}). The  kernel function $K_t$ is a Gaussian kernel with a variance equal to twice the weighted empirical variance of the set $S^{(t-1)}$ \cite{Beaumont2009}. The algorithm stops when the sample $S^{(T)}$ is generated i.e the target $\epsilon_T$ is reached. See Algorithm \ref{pmc} for details. 

\subsubsection*{Weights correcting the kernel sampling bias}

As pointed out by \cite{Beaumont2009}, the newly generated particles $\theta_i^{(t)}$ in a sequential procedure are no more drawn from the prior distribution but from a specific probability density $d_i^{(t)}$ that depends on the particles selected at the previous step and on the chosen kernel. This introduces a bias in the procedure. This bias should be corrected by attributing a weight equal to $\pi(\theta_i^{(t)}) / d_i^{(t)}$ to each newly generated particle $\theta_i^{(t)}$. 

The density of probability $d_i^{(t)}$ to generate particle $\theta_i^{(t)}$ at step $t$ is given by the sum of the probabilities to reach $\theta_i^{(t)}$ from one of the $N$ particles of the previous step times their respective weights:

\begin{equation}
	d_i^{(t)} \propto \sum_{j=1}^{N} w_{j}^{(t-1)} \sigma_{t-1}^{-1}\varphi \left(\sigma_{t-1}^{-1} (\theta_{i}^{(t)}-\theta_{j}^{(t-1)})\right)
\end{equation}

where $\varphi(x)=\frac{1}{\sqrt{2\pi}} e^{-\frac{x^2}{2}}$ is the kernel function.

This yields the expression of the weight $w_{i}^{(t)}$ to be attributed to the newly drawn particle $\theta_i^{(t)}$:
\begin{equation}
	\label{weightPMC}
	w_{i}^{(t)} \propto  \frac{\pi(\theta_{i}^{(t)})}{ \sum_{j=1}^{N} w_{j}^{(t-1)} \sigma_{t-1}^{-1}\varphi \left(\sigma_{t-1}^{-1} (\theta_{i}^{(t)}-\theta_{j}^{(t-1)})\right)}
\end{equation}

\subsubsection*{Limitations of the PMC algorithm}

The major problem in the PMC algorithm is to define the decreasing sequence of tolerance levels $\{\epsilon_1,...,\epsilon_T\}$ to get close to an optimal gain in computing time. If the decrease in tolerance values is too sharp or too shallow, the benefits of the importance sampling procedure has good chance to be lower than what could be possible. In the following, we will indeed demonstrate that our algorithm leads to a sequence of tolerance levels which clearly outperforms an arbitrary choice for the sequence of tolerance levels.

\subsection*{The RSMC and the SMC algorithms}

In \cite{Drovandi2011} and \cite{DelMoral2012} the authors proposed two methods to determine "on-line" the sequence of tolerance levels. The main idea is to define the $\epsilon_t$ value with the previous sample $S^{(t-1)}$. In the RSMC algorithm of \cite{Drovandi2011}, $\epsilon_t$ is defined as a quantile of the $\rho(x,y)$ values of the previous sample $S^{(t-1)}$ (see Algorithm \ref{rsmc} for details). In the SMC algorithm of \cite{DelMoral2012}, $\epsilon_t$ is computed so that the effective sample size of the particles is reduced by a constant factor at each time step (see Algorithm \ref{smc} for details). 

A second difference between the PMC and the RSMC/SMC algorithms concerns the proposal distribution. The RSMC and the SMC algorithms use a MCMC kernel to move the particles. At step $t$, a new particle $\theta_i^{(t)}$ is generated using a MCMC kernel $\theta_i^{(t)}\sim K_t(\theta|\theta*)$ where $\theta^*$ is randomly draw from $S^{(t-1)}$ with probability $(w_i^{(t-1)})_{i = 1,..,N}$. This weight $(w_i^{(t)})_{i = 1,..,N}$ is equal to 1 if the particle $\theta_i^{(t)}$ satisfies $\epsilon_t$, and to 0 otherwise. The jump is accepted with probability, $p_{acc}$, based on the Metropolis-Hastings ratio (Eq. \ref{ratio}). 

\begin{equation}
	\label{ratio}
1 \wedge \displaystyle{\frac{\pi(\theta_i^{(t)})K_t(\theta^*|\theta_i^{(t)})}{\pi(\theta^*)K_t(\theta_i^{(t)}|\theta^*)}\mathds{1}_{\rho(x,y)\leq \epsilon_{t}}}
\end{equation}

where $x \wedge y$ means the minimum of $x$ and $y$.

\subsubsection*{Limitations of the RSMC and the SMC algorithms}

The MCMC kernel used in \cite{Drovandi2011} and \cite{DelMoral2012} to sample new values $\theta_j^{(t)}$ has a significant drawback in our view: it can lead to particle duplications. Indeed, each time the MCMC jumps from a particle to a new one which is not accepted, the initial particle is kept in the new sample of particles. When this occurs several times with the same initial particle, this particle appears several times in the new sample. The number of such "duplicated" particles can grow and strongly deteriorate the quality of the posterior, as illustrated below. To solve this problem, \cite{Drovandi2011} proposed to perform $R$ MCMC jump trials instead of one. $R$ evolves during the course of the algorithm (Eq. \ref{R}) since its value is chosen such that there is a probability of $1-c$ that the particle gets moved at least once where $c=0.01$ in \cite{Drovandi2011}. To circumvent the problem of particle duplications \cite{DelMoral2012} proposed to resample the parameter values when too many are duplicated. In \cite{DelMoral2012} the authors also proposed to run the model $M$ times for each particle, in order to decrease the variance of the acceptance ratio of the MCMC jump. However, all these solutions increase the number of model runs, going against the initial benefit of using sequential samples.

\begin{equation}
	\label{R}
  R=\displaystyle{\frac{\log(c)}{\log(1-p_{acc})}}
\end{equation}

\section*{Adaptive population Monte-Carlo approximate Bayesian computation}

\subsection*{Overview of the APMC algorithm}

The APMC algorithm follows the main principles of the sequential ABC, and defines on-line the tolerance level at each step like in \cite{Wegmann2010}, \cite{Drovandi2011} and \cite{DelMoral2012}. For each tolerance level $\epsilon_t$, it generates a sample $S^{(t)}$ of particles and computes their associated weights. This weighted sample approximates the posterior distribution, with an increasing approximation quality as $\epsilon_t$ decreases. 
Suppose the APMC reached step $t-1$, with a sample $S^{(t-1)}$ of $N_\alpha=\lfloor \alpha N \rfloor$ particles and their associated weights $(\theta_i^{(t-1)},w_i^{(t-1)})_{i = 1,..,N_{\alpha}}$, the main features of the APMC are (see Algorithm \ref{apmc} for details): 
\begin{itemize}
	\item the algorithm generates $N-N_\alpha$ particles $(\theta_j^{(t-1)})_{j = N_{\alpha}+1,..,N}$ where $\theta_j^{(t-1)}\sim \mathcal{N}(\theta_j^*,\sigma_{(t-1)}^2)$, the seed $\theta_j^*$ is randomly drawn from the weighted set $(\theta_i^{(t-1)},w_i^{(t-1)})_{i = 1,..,N_{\alpha}}$ and the variance $\sigma_{(t-1)}^2$ of the Gaussian kernel $\mathcal{N}(\theta_j^*,\sigma_{(t-1)}^2)$  is twice the empirical variance of the weighted set $(\theta_i^{(t-1)},w_i^{(t-1)})_{i = 1,..,N_\alpha}$, following \cite{Beaumont2009}.
	\item the weights $w_j^{(t-1)}$ of the new particles $(\theta_j^{(t-1)})_{j = N_{\alpha}+1,..,N}$ are computed so that these new particles can be combined with the sample $S^{(t-1)}$ of the previous step without causing a bias in the posterior distribution. These weights are given by Eq. \ref{weight} (see below).
	\item the algorithm concatenates the $N_\alpha$ previous particles $(\theta_i^{(t-1)})_{i = 1,..,N_{\alpha}}$ with the $N-N_\alpha$ new particles $(\theta_j^{(t-1)})_{j = N_{\alpha}+1,..,N}$, together with their associated weights and distances to the data. This constitutes a new set noted $S_{temp}^{(t)}=(\theta_i^{(t)},w_i^{(t)},\rho_i^{(t)})_{i = 1,..,N}$.
	\item the next tolerance level $\epsilon_{t}$ is determined as the first $\alpha-$quantile of the $(\rho_i^{(t)})_{i = 1,..,N}$.
	\item the new sample $S^{(t)}=(\theta_i^{(t)},w_i^{(t)})_{i = 1,..,N_{\alpha}}$ is then constituted from the $N_{\alpha}$ particles of $S_{temp}^{(t)}$ satisfying the tolerance level $\epsilon_{t}$. 
	\item if the proportion $p_{acc}$ of particles satisfying the tolerance level $\epsilon_{t-1}$ among the $N-N_{\alpha}$ newly generated particles is below a chosen value $p_{acc_{min}}$, the algorithm stops, and its result is $(\theta_i^{(t)})_{i = 1,..,N_\alpha}$ with their associated weights.
	\end{itemize}

Note that in our algorithm, to get a number $N_\alpha$ of retained particles for the next step, the choice of $\epsilon_t$ is heavily constrained: it has to be at least equal to the first $\alpha-$quantile of the $(\rho_i^{(t)})_{i = 1,..,N}$ and smaller than the immediately superior $(\rho_i^{(t)})$ value. We chose to fix it to the first $\alpha-$quantile for simplicity. This choice also ensures that the tolerance level decreases from one iteration to the next: in the worst case where $p_{acc}=0$ (no newly simulated particles accepted), $\epsilon_t=\epsilon_{t-1}$. Our algorithm does not use a MCMC kernel and avoids duplicating particles. It requires a reweighting step in $O(N_\alpha^2)$ instead of $O(N_\alpha)$ in \cite{Drovandi2011}, but in our perspective, this computational cost is supposed negligible compared with the cost of running the model. 

\subsection*{Weights correcting the kernel sampling bias}

For the APMC algorithm the density of probability $d_i^{(t)}$ to generate particle $\theta_i^{(t)}$ at step $t$ is:

\begin{equation}
	d_i^{(t)}= \sum_{j=1}^{N_\alpha} \frac{w_{j}^{(t-1)}}{\sum_{k=1}^{N_ \alpha} w_{k}^{(t-1)}} \sigma_{t-1}^{-1}\varphi \left(\sigma_{t-1}^{-1} (\theta_{i}^{(t)}-\theta_{j}^{(t-1)})\right)
\end{equation}

where $\varphi(x)=\frac{1}{\sqrt{2\pi}} e^{-\frac{x^2}{2}}$ is the kernel function.

This yields the expression of the weight $w_{i}^{(t)}$ to be attributed to the newly drawn particle $\theta_i^{(t)}$:

\begin{equation}
	\label{weight}
	w_{i}^{(t)}=  \frac{\pi(\theta_{i}^{(t)})}{ \sum_{j=1}^{N_\alpha} \left(w_{j}^{(t-1)}/\sum_{k=1}^{N_ \alpha} w_{k}^{(t-1)}\right) \sigma_{t-1}^{-1}\varphi \left(\sigma_{t-1}^{-1} (\theta_{i}^{(t)}-\theta_{j}^{(t-1)})\right)}
\end{equation}

This formula differs from the scheme of \cite{Beaumont2009} where the weights need only to be proportional to Eq. \ref{weight} at each step. Since we want to concatenate particles obtained at different steps of the algorithm (while \cite{Beaumont2009} generate the sample at step $t$ from scratch), we need the scaling of weights to be consistent across the different steps of the algorithm. Using the weight of Eq. \ref{weight} guarantees the correction of the sampling bias throughout the APMC procedure and ensures that the $N_\alpha$ weighted particles $\theta_i^{(t)}$ produced at the $t$-th iteration follow the posterior distribution $\pi\left(\theta|\rho(x,y)<\epsilon_t\right)$.

\begin{figure*}
	\centering
		\includegraphics[width=\linewidth]{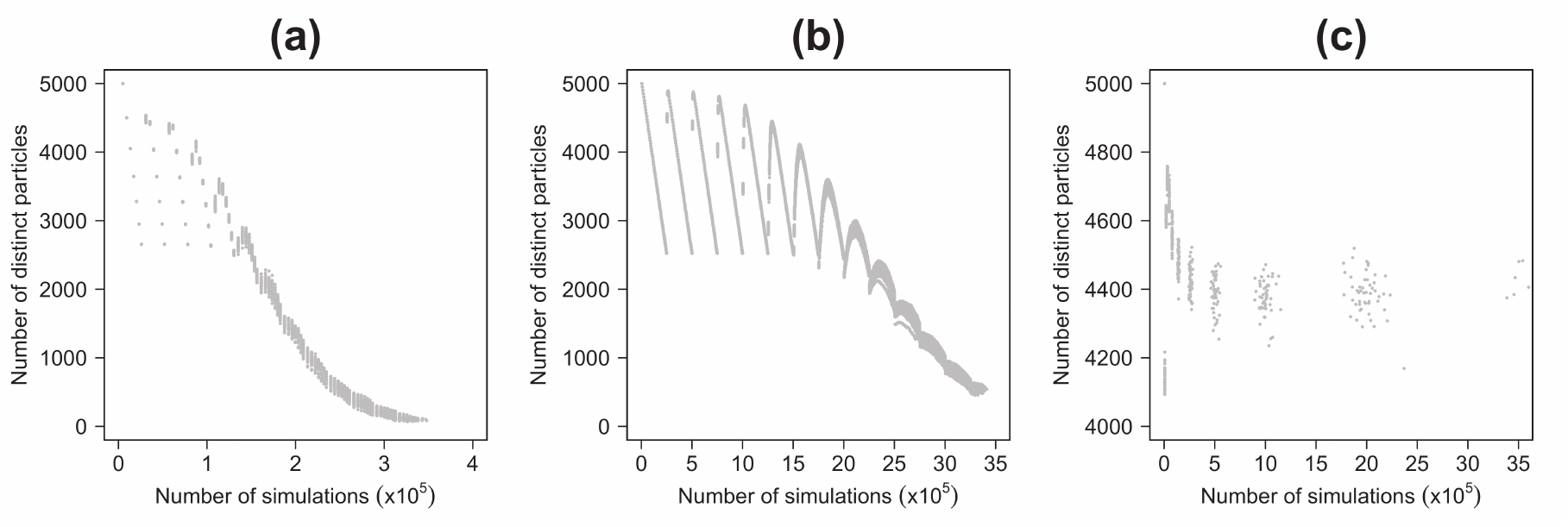}	
	\caption{Number of distinct particles in a sample of $N=5000$ particles during the course of the SMC and RSMC algorithms applied to the toy example; In all three panels we plot the results obtained for 50 executions of the algorithm. (a) SMC with 
		$\alpha=0.9$ and $M=1$; (b) SMC with $\alpha=0.99$ and $M=1$; (c) RSMC with $\alpha=0.5$. In all three panels, the tolerance target is equal to $0.001$. }
	\label{Fig1}
\end{figure*}

\subsection*{The stopping criterion}

We stop the algorithm when the proportion of "accepted" particles (Eq. \ref{pacc}) among the $N-N_{\alpha}$ new particles is below a predetermined threshold $p_{acc_{min}}$. This choice of stopping rule ensures that additional simulations would only marginally change the posterior distribution. Note that this stopping criterion will be achieved even if $p_{acc_{min}}=0$, this ensures that the algorithm converges. We present a formal proof of this assertion in Appendix B.

\begin{equation}
p_{acc}(t)=\frac{1}{N-N_\alpha}\sum_{k=N_\alpha+1}^N \mathds{1}_{\rho_{k}^{(t-1)} < \epsilon_{t-1}} \label{pacc}
\end{equation} 

\section*{Experiments on a toy example}
\label{ToyEx}

We consider four algorithms: APMC, PMC, the SMC and the RSMC. Their implementations in R \cite{R} are available \footnote[1]{\url{http://motive.cemagref.fr/people/maxime.lenormand/script_r_toyex}}. We compare them on the toy example studied in \cite{Sisson2007} where $\pi(\theta) = \mathcal{U}_{[-10,10]}$ and $f(x|\theta) \sim \frac{1}{2}\phi\left(\theta,\frac{1}{100}\right)+\frac{1}{2}\phi\left(\theta,1\right)$ where $\phi\left(\mu,\sigma^2\right)$ is the normal density of mean $\mu$ and variance $\sigma^2$. In this example, we consider that $y=0$ is observed, so that the posterior density of interest is proportional to
$\left(\phi\left(0,\frac{1}{100}\right)+\phi\left(0,1\right)\right)\pi(\theta)$.

We structure the comparisons on two indicators: the number of simulations performed during the application of the algorithms, and the $\mathbb{L}_2$ distance between the exact posterior density and the histogram of particle values obtained with the algorithms. This $\mathbb{L}_2$ distance is computed on the 300-tuple obtained by dividing the support $[-10,10]$ into 300 equally-sized bins. We choose the $\mathbb{L}_2$ distance to compare the sample to the true posterior because it is a well-known accuracy measure easy to compute and a good indicator to compare different methods.

We choose $N=5000$ particles and a target tolerance level equal to $0.01$. For the PMC algorithm we use a decreasing sequence of tolerance levels from $\epsilon_1=2$ down to $\epsilon_{11}=0.01$. For the SMC algorithm, we use $3$ different values for $\alpha$: $\{0.9,0.95,0.99\}$ and $M=1$ as in \cite{DelMoral2012}. For the RSMC algorithm we use $\alpha=0.5$ as in \cite{Drovandi2011}. To explore our algorithm, we test $9$ different values for $\alpha$: $\{0.1,0.2,0.3,0.4,0.5,0.6,0.7,0.8,0.9\}$, and $4$ different values for $p_{acc_{min}}$: $\{0.01,0.05,0.1,0.2\}$. In each case, we perform $50$ times the algorithm, and compute the average and standard deviation of the two indicators: the total number of simulations and the $\mathbb{L}_2$ distance between the exact posterior density and the histogram of particle values. We used as kernel transition a normal distribution parameterized with twice the weighted variance of the previous sample, as in \cite{Beaumont2009}.

We report below the effects of varying $\alpha$ and $p_{acc_{min}}$ on the performance of our algorithm, and compare it with the PMC, SMC and RSMC algorithms.

\subsection*{Particle duplication in SMC and RSMC}

\begin{figure*}
	\centering
	\includegraphics[scale=0.55]{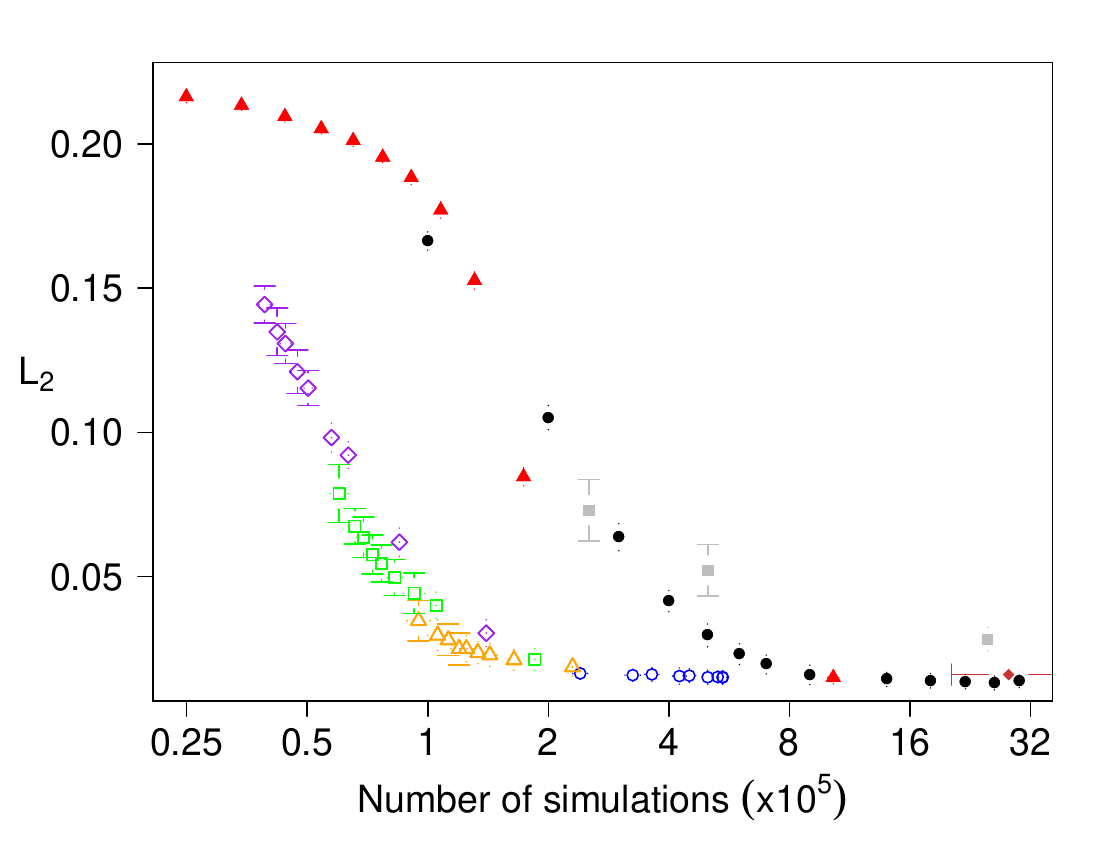}
	\caption{Posterior quality ($\mathbb{L}_2$) versus computing cost (number of simulations) averaged over 50 replicates. Vertical and horizontal bars 	
	represent the standard deviations among replicates. Algorithm parameters used for APMC: $\alpha$ in $\{0.1,0.2,0.3,0.4,0.5,0.6,0.7,0.8,0.9\}$ and 
	$p_{acc_{min}}$ in $\{0.01,0.05,0.1,0.2\}$. Blue circles are used for $p_{acc_{min}}=0.01$, orange triangles for $p_{acc_{min}}=0.05$, green squares 
	for $p_{acc_{min}}=0.1$, and purple diamonds for $p_{acc_{min}}=0.2$. PMC: red plain triangles for a sequence of tolerance levels from $\epsilon_1 = 
	2$ down to $\epsilon_{11} = 0.01$. SMC: grey plain square for $\alpha$ in $\{0.9,0.95,0.99\}$ (from left to right), $M=1$ and a $\epsilon$ target equal to 0.01. RSMC: 
	brown plain diamond for $\alpha=0.5$ and a $\epsilon$ target equal to 0.01. Results obtained with a standard rejection-based ABC algorithm are 
	depicted with black plain circles.}
	\label{Fig2}
\end{figure*}

The number of distinct particles decreases during the course of the SMC algorithm whatever the value of $\alpha$, as shown on Fig. \ref{Fig1}a-b. The oscillations of the number of distinct particles are caused by the resampling step in the SMC algorithm (see \cite{DelMoral2012}), but they are not sufficient to counterbalance the overall decrease. This decrease deteriorates the posterior approximation as shown on Fig. \ref{Fig2}. For the RSMC algorithm, the initial oscillation of the number of particles is due to the initial value of $R$, initially set to 1, but which quickly evolves towards a value ensuring a relatively constant number of distinct particles. This number of distinct particles is maintained at a reasonably high level (Fig. \ref{Fig1}c), but this has a cost in terms of the number of required model runs (see Fig. \ref{Fig2}). Note that the APMC and the PMC algorithms keep $N$ distinct particles.

	\begin{figure*}
		\begin{minipage}{.5\linewidth}
		  \centering
		  \includegraphics[scale=0.31]{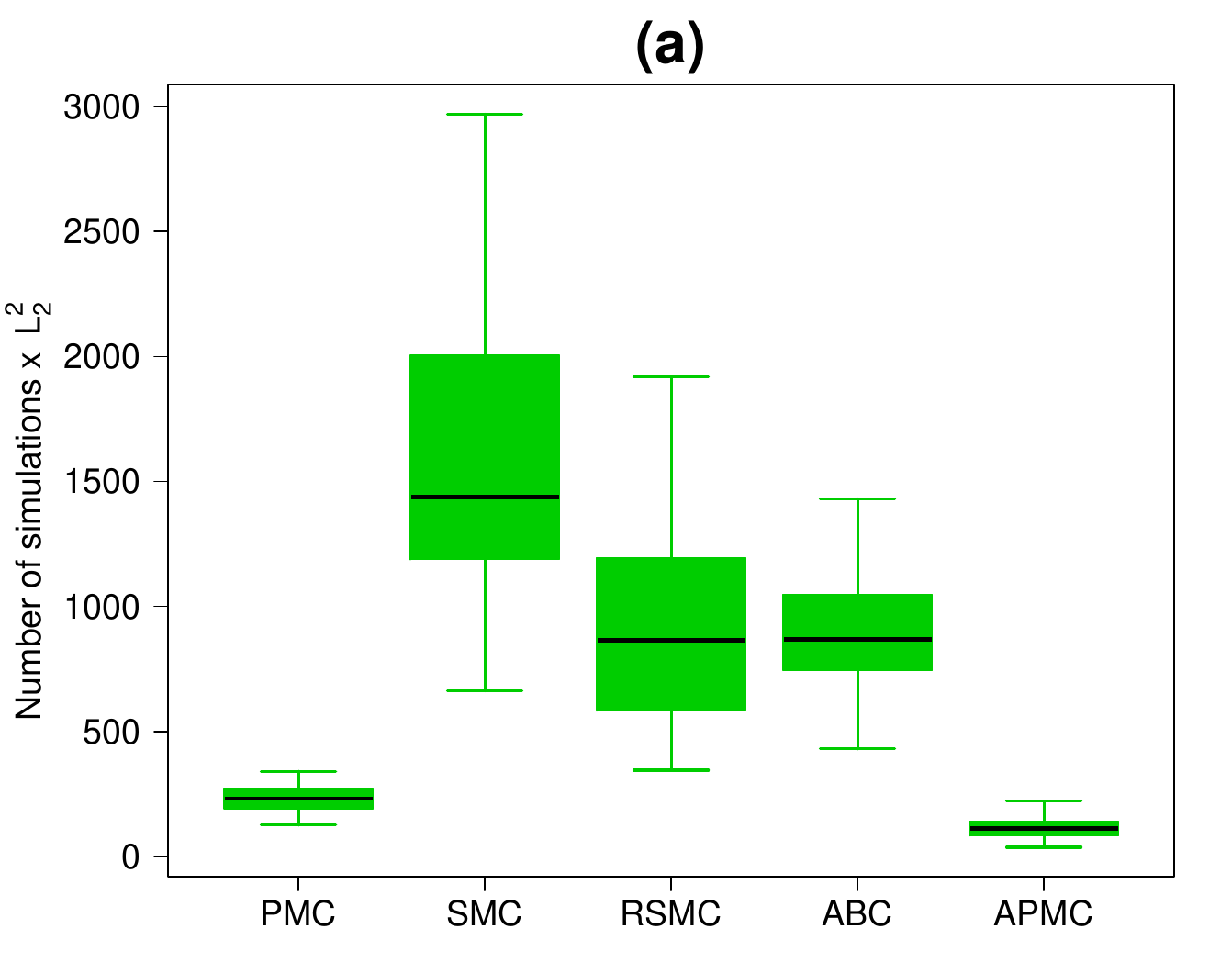}
		\end{minipage}\hfill
		\begin{minipage}[r]{.5\linewidth}
		   \centering
			\includegraphics[scale=0.31]{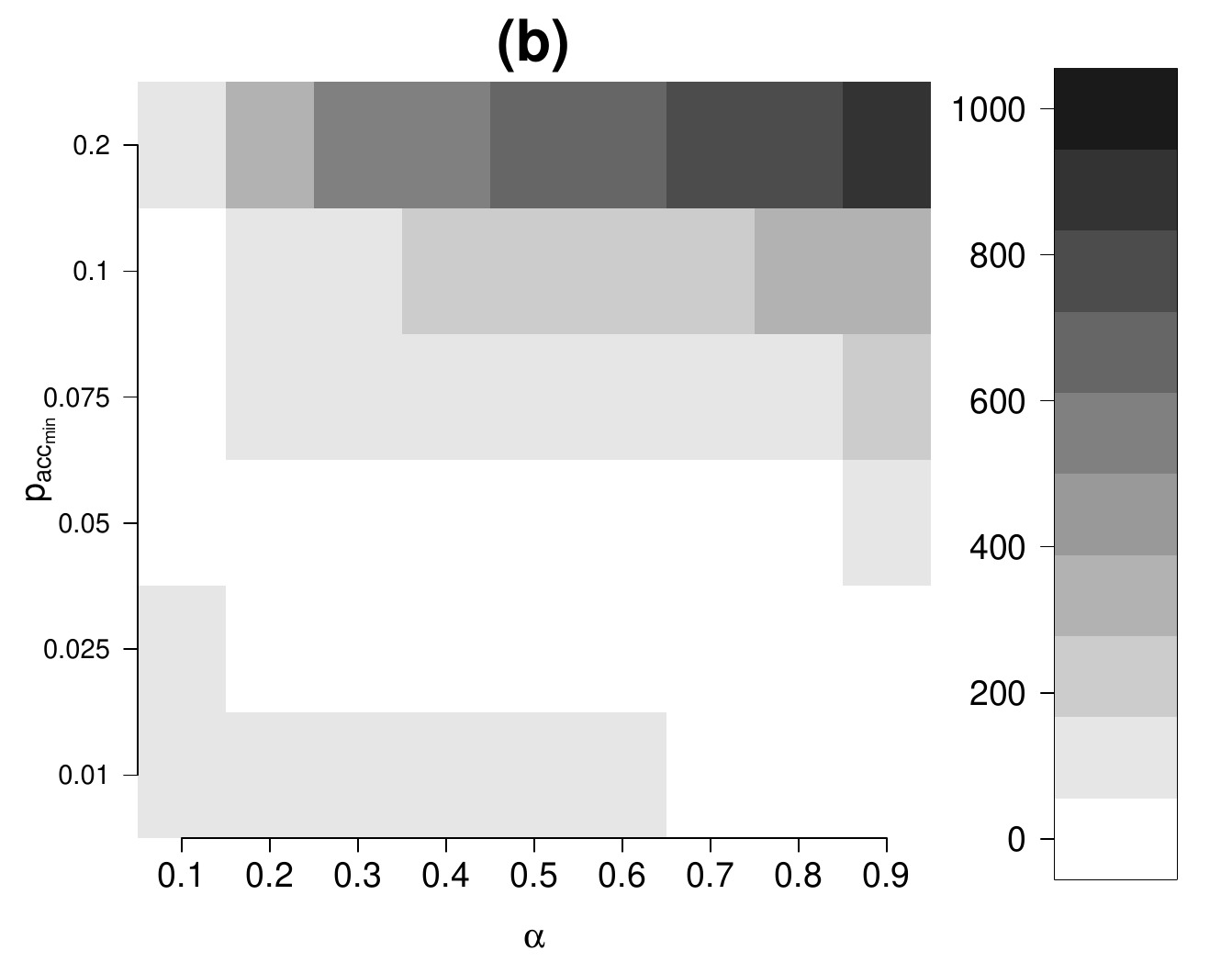}
		\end{minipage}		
		\caption{(a) Boxplot of the criterion ``squared $\mathbb{L}_2$ distance times the number of simulations'' for the different ABC algorithms. APMC: for $\alpha$ in $\{0.1,0.2,0.3,0.4,0.5,0.6,0.7,0.8,0.9\}$ and $p_{acc_{min}}=0.01$; SMC: for $\alpha$ in 
		$\{0.9,0.95,0.99\}$, $M=1$ and a $\epsilon$ target equal to 0.01; RSMC: for $\alpha=0.5$ and a $\epsilon$ target equal to 0.01; ABC: for a 
		$\epsilon$ target equal to 0.01; PMC: for a sequence of tolerance levels from $\epsilon_1 = 2$ to $\epsilon_{11} = 0.01$. (b) Criterion ``squared $\mathbb{L}_2$ distance times the number of simulations'' in the APMC algorithm for the different values of $\alpha$ and $p_{acc_{min}}$. Each cell depicts the average of the criterion over the 50 performed replicates of the APMC.}
		\label{Fig3}
	\end{figure*}

\vspace{5cm}

\begin{table*}
\caption{SimVillages parameter descriptions}
\label{param}
\begin{center}
\begin{tabular}{ccc} 
\hline\noalign{\smallskip}
\textbf{Parameters} & \textbf{Description} & \textbf{Range}\\
\noalign{\smallskip}
\hline
\noalign{\smallskip}
$\theta_1$ & Average number of children per woman& $[0,4]$\\
$\theta_2$ & Probability to accept a new residence for a household & $[0,1]$\\
$\theta_3$ & Probability to make couple for two individuals & $[0,1]$\\
$\theta_4$ & Probability to split for a couple in a year & $[0,0.5]$\\
\hline
\end{tabular}
\end{center}
\end{table*}

\begin{table*}
\caption{Summary statistic descriptions}
\label{sumstat}
\begin{center}
\begin{tabular}{ccc}
\hline\noalign{\smallskip}
\textbf{Summary statistic} & \textbf{Description} & \textbf{Measure of discrepancy}\\
\noalign{\smallskip}
\hline
\noalign{\smallskip}
$S_1$ & Number of inhabitants in 1999 & $\mathbb{L}_1$ distance\\
$S_2$ & Age distribution in 1999 & $\chi^2$ distance\\
$S_3$ & Household type distribution in 1999 & $\chi^2$ distance\\
$S_4$ & Net migration in 1999 & $\mathbb{L}_1$ distance\\
$S_5$ & Number of inhabitants in 2006 & $\mathbb{L}_1$ distance\\
$S_6$ & Age distribution in 2006 & $\chi^2$ distance\\
$S_7$ & Household type distribution in 2006 & $\chi^2$ distance\\
$S_8$ & Net migration in 2006 & $\mathbb{L}_1$ distance\\
\hline
\end{tabular}
\end{center}
\end{table*}

\subsection*{Influence of parameters on APMC}

The values of $\alpha$ and $p_{acc_{min}}$ have an impact on the studied indicators. We find that smaller $\alpha$ and $p_{acc_{min}}$ improve the quality of the approximation (smaller $\mathbb{L}_2$ distance), and increase the total number of model runs, with $p_{acc_{min}}$ having the largest effect (Fig. \ref{Fig2}). With a large $\alpha$, the tolerance levels decrease slowly and there are numerous steps before the algorithm stops. In this toy example, our simulations show that all explored sets of ($\alpha$ , $p_{acc_{min}}$) such that $p_{acc_{min}} < 0.1$ give good results for the criterion $\mbox{\textit{Number of simulations} } \times \mbox{ } \mathbb{L}_2^2$ (Fig. \ref{Fig3}b). Large $\alpha$ provide slightly better results for small $p_{acc_{min}}$ while small $\alpha$ provide slightly better results for large $p_{acc_{min}}$ (Fig. \ref{Fig3}b). On this toy example it appears that intermediate values of $\alpha$ and $p_{acc_{min}}$ ($0.3 \leq \alpha \leq 0.7$ and $0.01 \leq p_{acc_{min}} \leq 0.05$), present a good compromise between number of model runs and the quality of the posterior approximation.

\subsection*{Comparing performances}

Whatever the value of $\alpha$ and $p_{acc_{min}}$, the APMC algorithm always yields better results than the other three algorithms. It requires between $2$ and $8$ times less simulations to reach a given posterior quality $\mathbb{L}_2$ (Fig. \ref{Fig2}). Furthermore, good approximate posterior distributions are very quickly obtained (Fig. \ref{Fig2}). The compromise between simulation speed and convergence level can also be illustrated using the criterion $\mbox{\textit{Number of simulations} } \times \mbox{ } \mathbb{L}_2^2$ \cite{Glynn1992}. This criterion is smaller for the APMC algorithm (Fig. \ref{Fig3}a). 

\section*{Application to the model SimVillages}

In this section, we check if our algorithm still performs better than the PMC, the RSMC and the SMC when applied to an individual-based social model developed during the European project PRIMA\footnote[2]{PRototypical policy Impacts on Multifunctional Activities in rural municipalities - EU 7th Framework Research Programme; 2008-2011; \url{https://prima.cemagref.fr/the-project}}. The aim of the model is to simulate the effect of a scenario of job creation (or destruction) on the evolution of the population and activities in a network of municipalities.

\begin{figure*}
	\centering
	\includegraphics[scale=0.48]{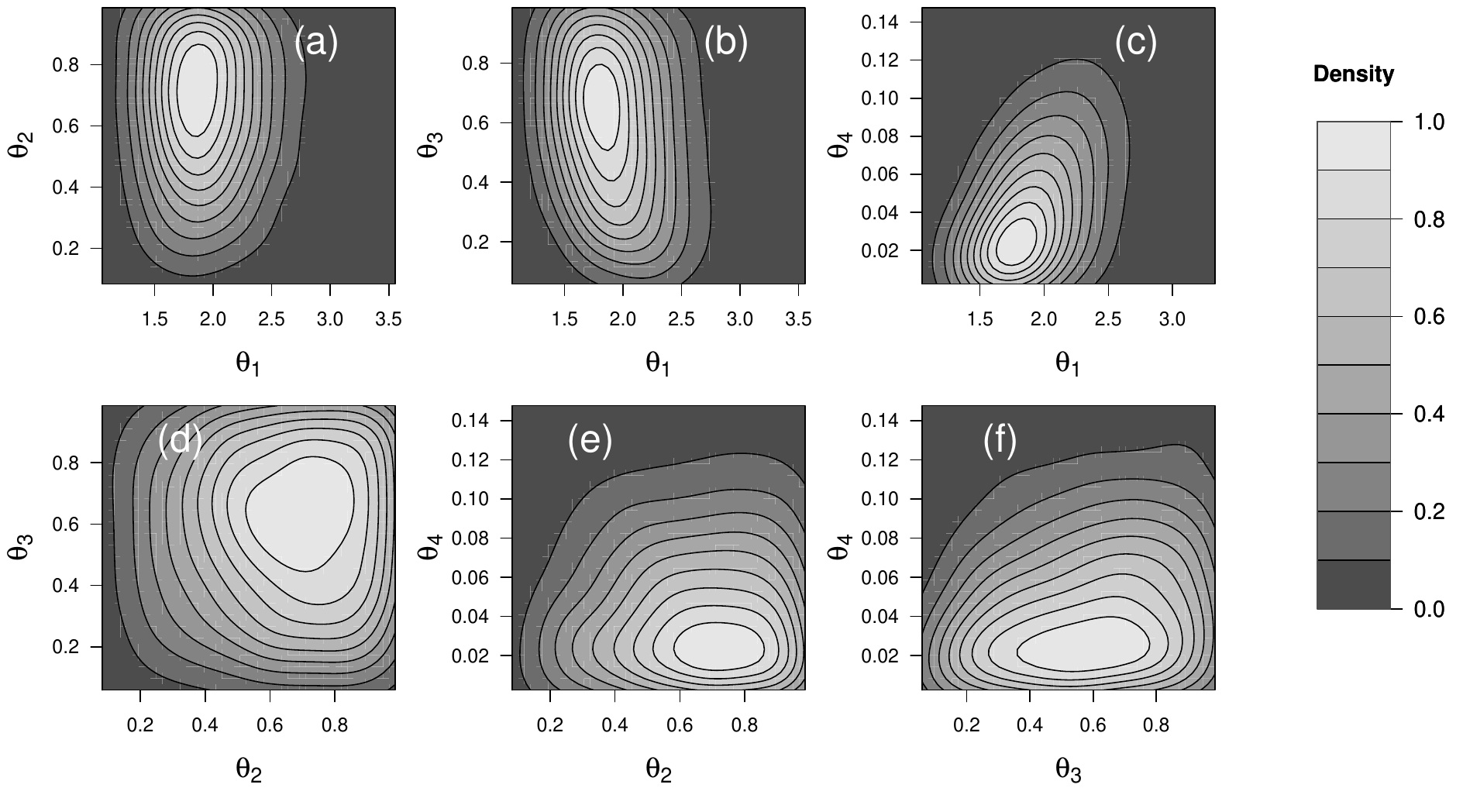} 
	\caption{Contour plot of the bivariate joint densities of $\theta_i$ and $\theta_j$ obtained with our algorithm, and with $\alpha=0.5$ and 		  	   
	$p_{acc_{min}}=0.01$; (a) $\theta_1$ and $\theta_2$; (b) $\theta_1$ and $\theta_3$; (c) $\theta_1$ and $\theta_4$; (d) $\theta_2$ and $\theta_3$; (e) 
	$\theta_2$ and $\theta_4$; (f) $\theta_3$ and $\theta_4$.}
  \label{Fig4}
\end{figure*}

\subsection*{Model and data} 

The model simulates the dynamics of virtual individuals living in 7 interconnected villages in a rural area of Auvergne (a region of Central France). A single run of the model SimVillages with seven rural municipalities takes about $1.4$ seconds on a desktop machine (PC Intel 2.83 GHz). The dynamics include demographic change (aging, marriage, divorce, births and deaths), activity change (change of jobs, unemployment, inactivity, retirement), and movings from one municipality to another or outside of the set. The model also includes a dynamics of creation / destruction of jobs of proximity services, derived from the size of the local population. More details on the model can be found in \cite{Huet2011}. The individuals (about 3000) are initially generated using the 1990 census data of the National Institute of Statistics and Economic Studies ($INSEE$), some of them are given a job type and a location for this job (in a municipality of the set or outside), they are organised in households living in a municipality of the set. The model dynamics is mostly data driven, but four parameters cannot be directly derived from the available data. They are noted $\theta_p$ for $1\leq p \leq 4$, described in Table \ref{param}.

We use our algorithm to identify the distribution of the four parameters for which the simulations, initialized with the 1990 census data, satisfy matching criteria with the data of the 1999 and 2006 census. The set of summary statistics $\{S_m\}_{1 \leq m \leq M}$ and the associated discrepancy measure used $\rho_m$ are described in Table \ref{sumstat}. We note $S_m$ the simulated summary statistics and $S_m^{'}$ the observed statistics. The eight summary statistics are normalized (variance equalization) and they are combined using the infinity norm (Eq. \ref{norminf}):
\begin{equation}
{\|{(\rho_m(S_m,S_m^{'}))}_{1 \leq m \leq M}\|}_{\infty} = \sup_{1 \leq m \leq M} \rho_m(S_m,S_m^{'})
\label{norminf}
\end{equation}
We first generate a sample of length $N$ from the prior $\mathcal{U}_{[a,b]}$, where $[a,b]$ is available for each parameter in Table \ref{param}, with a Latin hypercube \cite{Carnell2009} and we select the best $N_\alpha$ particles. To move the particles, we use as kernel transition a multivariate normal distribution parameterized with twice the weighted variance-covariance matrix of the previous sample \cite{Filippi2011}.

As in the section \ref{ToyEx}, we perform a parameter study and compare APMC with its three competitors. For APMC, $\alpha$ varies in ($\{0.3,0.5,0.7\}$) and $p_{acc_{min}}$ in ($\{0.01,0.05,0.1,0.2\}$), and we set $N_\alpha=5000$ particles. For the PMC, SMC and RSMC we also set $N=5000$ particles and a tolerance level target equal to $1.4$. The tolerance value $\epsilon=1.4$ corresponds to the average final tolerance value we obtain with APMC for $p_{acc_{min}}=0.01$. Note that otherwise this final tolerance is difficult to set properly and a worse choice for this value would have lead to worse performances of these algorithms. For the PMC algorithm, we use the decreasing sequence of tolerance levels $\{3,2.5,2,1.7,1.4\}$. For the SMC algorithm, we use $3$ different values for the couple $(\alpha,M)$: $\{(0.9,1),(0.99,1)$ $,(0.9,15)\}$. For the RSMC algorithm we use $\alpha=0.5$, as in \cite{Drovandi2011}. For each algorithm and parameter setting, we perform $5$ replicates.

We approximated posterior density (unknown in this case) with the original rejection-based ABC algorithm, starting with $N = 10,000,000$, selecting $7890$ particles below the tolerance level $\epsilon = 1.4$. 
 
To compute the $\mathbb{L}_2$ distance between posterior densities, we divided each parameter support into 4 equally sized bins, leading to a grid of $4^4=256$ cells, and we computed on this grid the sum of the squared differences between histogram values.

\subsection*{Study of APMC result}
APMC yields a unimodal approximate posterior distribution for the model SimVillages (Fig. \ref{Fig4}). Interestingly, parameters $\theta_1$ and $\theta_4$ are slightly correlated (Fig. \ref{Fig4}c). This is logical since they have contradictory effects on the number of children in the population. What is less straightforward is that we are able to partly tease apart these two effects with the available census data, since we get a peak in the approximate posterior distribution instead of a ridge.

\subsection*{Influence of parameters on APMC}

As for the toy example, we find that the intermediate values of $(\alpha,p_{acc_{min}})$ that we used lead to similar results (Fig. \ref{Fig5}c). In practice, we therefore recommend to use $\alpha=0.5$ and $p_{acc_{min}}$ between $0.01$ and $0.05$ depending on the wished level of convergence.

\subsection*{Comparing performances}

\begin{figure*}
		  \centering
		  \includegraphics[scale=0.3]{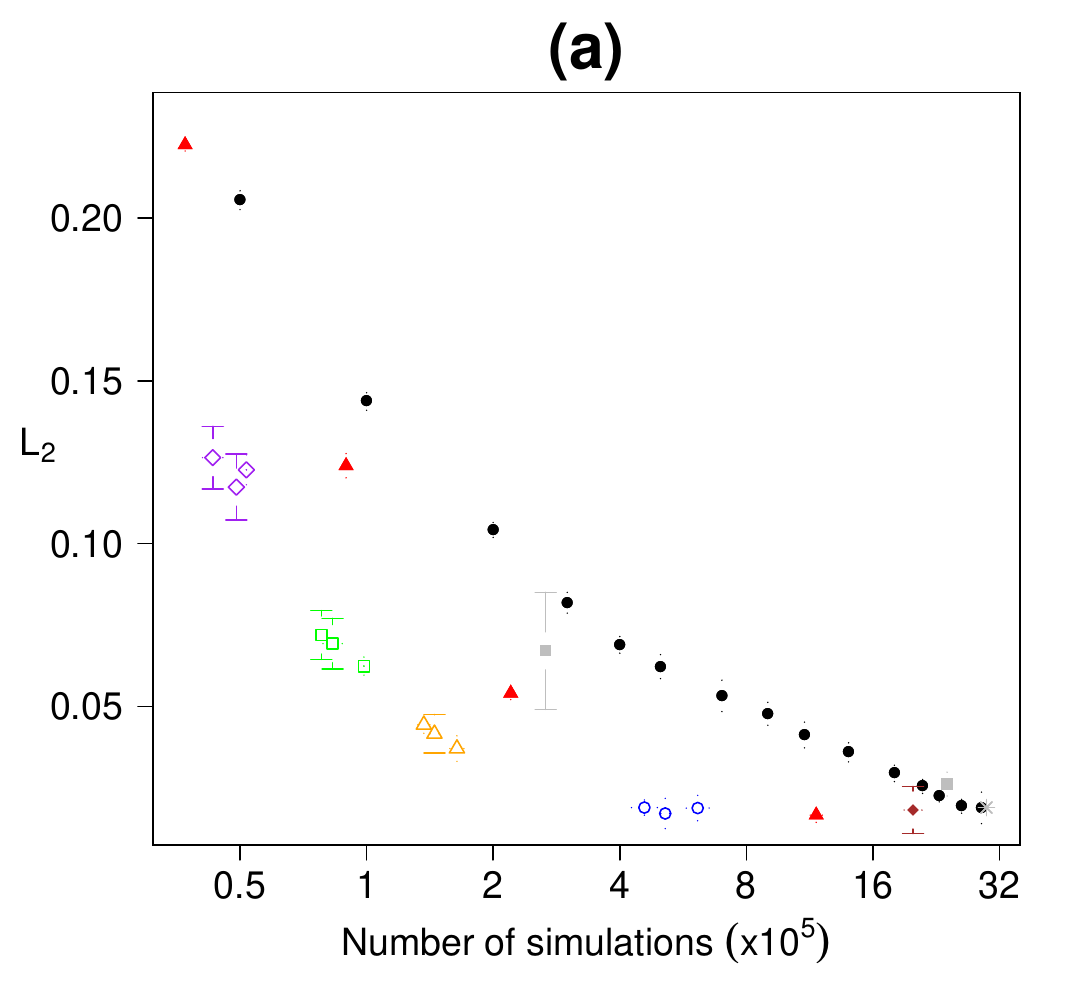}
			\includegraphics[scale=0.3]{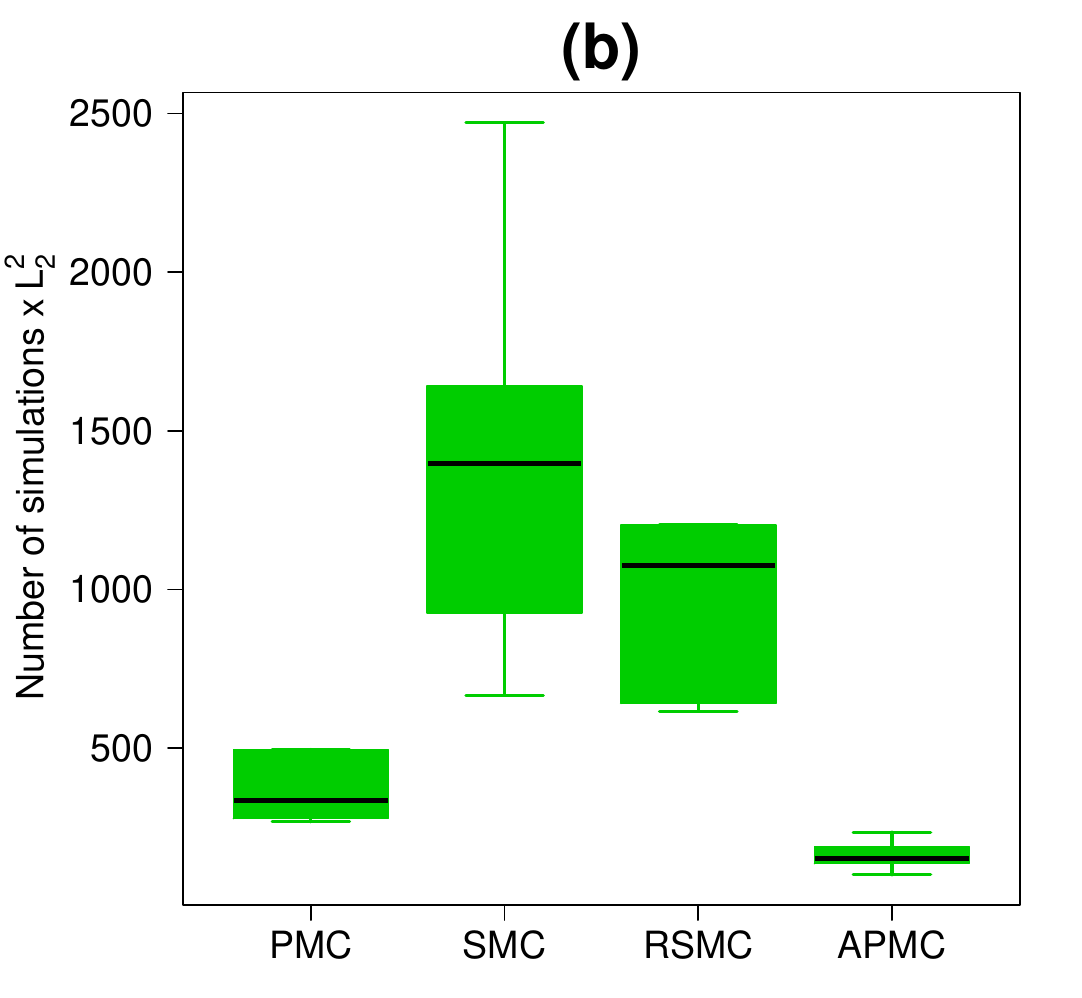}
			\includegraphics[scale=0.3]{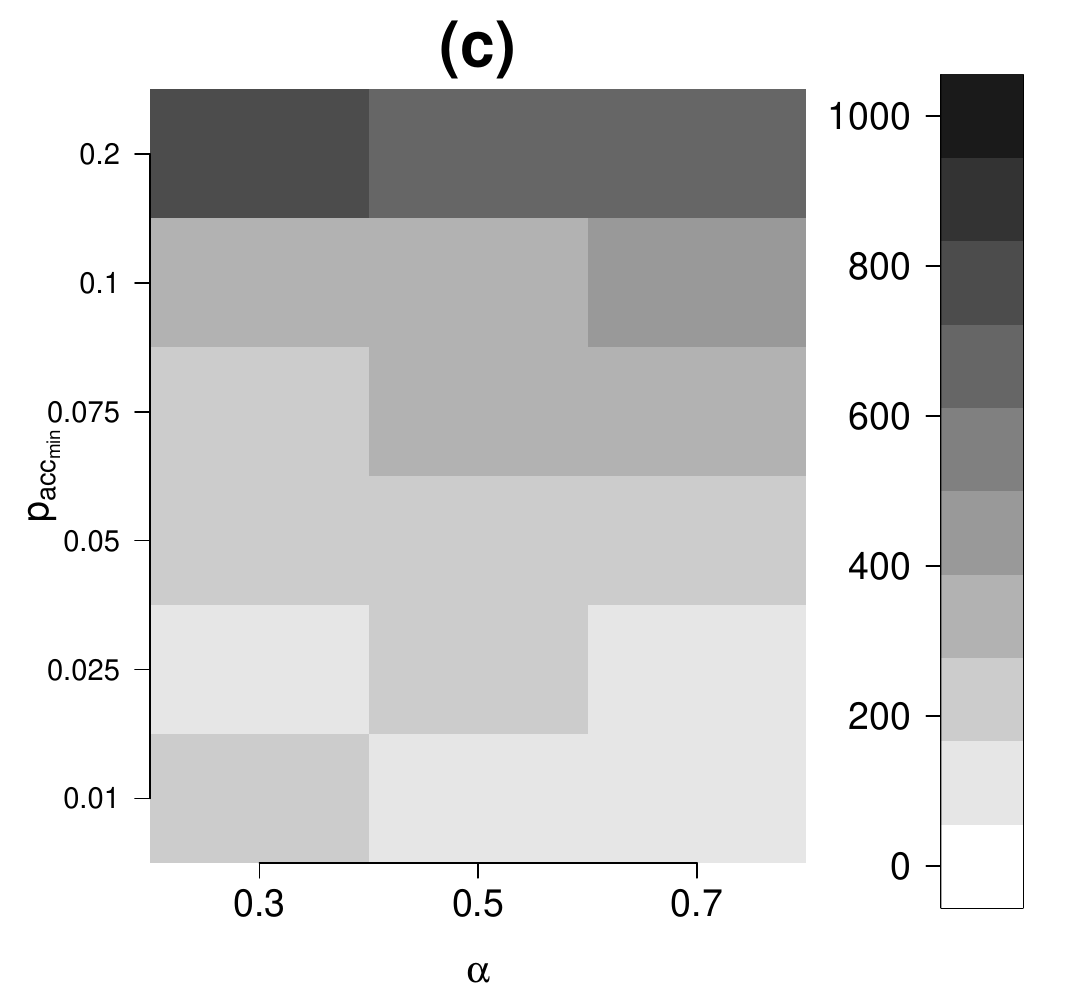}
		\caption{(a) Posterior quality ($\mathbb{L}_2$) versus computing cost (number of simulations) averaged over 5 replicates. Vertical and horizontal 	
		bars represent the standard deviations among replicates. Algorithm parameters used for APMC: $\alpha$ in $\{0.3,0.5,0.7\}$ and $p_{acc_{min}}$ in 
		$\{0.01,0.05,0.1,0.2\}$. Blue circles are used for $p_{acc_{min}}=0.01$, orange triangles for $p_{acc_{min}}=0.05$, green squares for 
		$p_{acc_{min}}=0.1$, and purple diamonds for $p_{acc_{min}}=0.2$. PMC: red plain triangles for a sequence of tolerance levels from $\epsilon_1 = 3$ 
	  to $\epsilon_{5} = 1.4$. SMC: grey plain square for $(\alpha,M)$ in $\{(0.9,1),(0.99,1)\}$, grey star for $(\alpha,M)=(0.9,15)$ 
	  and a $\epsilon$ target equal to 1.4. RSMC: brown plain diamond for $\alpha=0.5$ and a $\epsilon$ target equal to 1.4. Results obtained with a 
	  standard rejection-based ABC algorithm are depicted with black plain circles. (b) Boxplot of the criterion ``squared $\mathbb{L}_2$ distance times 
	  the number of simulations'' for the different algorithms. APMC: for $\alpha$ in $\{0.3,0.5,0.7\}$ and $p_{acc_{min}}=0.01$; SMC: for $(\alpha,M)$ 
	  in $\{(0.9,1),(0.99,1),(0.9,15)\}$ and a $\epsilon$ target equal to 0.01; RSMC: for $\alpha=0.5$ and a $\epsilon$ target equal to 
	  0.01; ABC: for a $\epsilon$ target equal to 1.4; PMC: for a sequence of tolerance levels from $\epsilon_1 = 3$ to $\epsilon_{5} = 1.4$. (c) 
	  Criterion ``squared $\mathbb{L}_2$ distance times the number of simulations'' in the APMC algorithm for the different values of $\alpha$ and 
	  $p_{acc_{min}}$. Each cell depicts the average of the criterion over the 5 performed replicates of the APMC.}
	  \label{Fig5}
\end{figure*}

APMC requires between 2 and 7 times less simulations to reach a given posterior quality than the other algorithms $\mathbb{L}_2$ (Fig. \ref{Fig5}a). Again, the gain in simulation number is progressive during the course of the algorithm. The $\mbox{\textit{Number of}}$  $\mbox{\textit{simulations} } \times \mbox{ } \mathbb{L}_2^2$ criterion is again smaller for the APMC algorithm (Fig. \ref{Fig5}b).

\section*{Discussion}

The good performances of APMC should of course be confirmed on other examples. Nevertheless we argue that they are due to the main assets of our approach:
\begin{itemize}
	\item We choose an appropriate reweighting process instead of a MCMC kernel, which corrects the sampling bias without duplicating particles;
	\item We define an easy to interpret stopping criterion that automatically defines the number of sequential steps.
\end{itemize}
Therefore, we can have some confidence in the good performances of APMC on other examples.

In the future, it would be interesting to evaluate this algorithm on models involving a larger number of parameters and/or multi-modal posterior distributions. Moreover, APMC could benefit from other improvements, in particular by performing a semi-automatic selection of informative summary statistics after the first ABC step \cite{Joyce2008,Fearnhead2011} and by using local regressions for post-processing the final posterior distribution \cite{Beaumont2002,Blum2010}. We did not perform such combinations in the present contribution, so that our algorithm is directly comparable with the three other sequential algorithms we looked at. However, they would be straightforward, because the different improvements concern different steps of the ABC procedure.

\vspace{0.5cm}
\section*{Acknowledgements}

This publication has been funded by the Prototypical policy impacts on multifunctional activities in rural municipalities collaborative project, European Union 7th Framework Programme (ENV 2007-1), contract no. 212345. The work of the first author has been funded by the Auvergne region.

\bibliographystyle{unsrt}  
\bibliography{APMC}

\onecolumngrid

\newpage
\section*{Appendix A: Description of the algorithms}

\begin{algorithm}
  {\hrulefill}
  \vspace*{-0.4cm}
	\caption{Likelihood-free rejection sampler (ABC)}
	{\vspace*{-0.25cm}\hrulefill}
	\label{ABC}
	\begin{algorithmic}
	 \STATE Given $N$ the number of particles
		\FOR {$i=1$ to $N$} 
			\REPEAT
				 \STATE Generate $\theta^* \sim \pi(\theta)$
				 \STATE Simulate $x \sim f(x|\theta^*)$
			\UNTIL $\rho(S(x),S(y))<\epsilon$	  
			\STATE Set $\theta_i=\theta^*$
		\ENDFOR
	\end{algorithmic}
	{\vspace*{-0.2cm}\hrulefill}	
\end{algorithm}

\begin{algorithm}
  {\hrulefill}
  \vspace*{-0.4cm}
	\caption{Population Monte Carlo Approximate Bayesian Computation (PMC)}
	{\vspace*{-0.25cm}\hrulefill}	
	\label{pmc}
	\begin{algorithmic}
	  \STATE Given $N$ the number of particles and a decreasing sequence of tolerance level $\epsilon_1 \geq...\geq \epsilon_T$,
	  \STATE For $t=1$,
			\FOR {$i=1$ à $N$} 
		  	\REPEAT 
					\STATE $\mbox{Simulate } \theta_{i}^{(1)} \sim \pi(\theta) \mbox{ and }x\sim f(x|\theta_{i}^{(1)})$ 
				\UNTIL $\rho(S(x),S(y))<\epsilon_1$
				\STATE $\mbox{Set } \displaystyle{w_{i}^{(1)}=\frac{1}{N}}$
			\ENDFOR
			\STATE $\mbox{Take } \sigma_2^2$ as twice the weighted empirical variance of $(\theta_{i}^{(1)})_{1\leq i \leq N}$
		\FOR {$t=2$ to $T$} 
			\FOR {$i=1$ to $N$} 
		  	\REPEAT 
					\STATE Sample $\theta_i^* \mbox{ from } \theta_{j}^{(t-1)} \mbox{ with probabilities } w_{j}^{(t-1)}$
					\STATE Generate $\theta_{i}^{(t)}| \theta_i^* \sim \mathcal{N}(\theta_i^*,\sigma_t^2) \mbox{ and } x\sim f(x|\theta_{i}^{(t)})$ 
				\UNTIL $\rho(S(x),S(y))<\epsilon_t$
				\STATE Set $\displaystyle{w_{i}^{(t)}\propto \frac{\pi(\theta_{i}^{(t)})}{\sum_{j=1}^N w_{j}^{(t-1)} \sigma_t^{-1}\varphi(\sigma_t^{-1} (\theta_{i}^{(t)}-\theta_{j}^{(t-1)}))}}$
			\ENDFOR
			\STATE Take $\sigma_{t+1}^2$ as twice the weighted empirical variance of $(\theta_{i}^{(t)})_{1\leq i \leq N}$
		\ENDFOR
		\STATE Where $\varphi(x)=\displaystyle{\frac{1}{\sqrt{2\pi}} e^{-\frac{x^2}{2}}}$
	\end{algorithmic}
	{\vspace*{-0.2cm}\hrulefill}		
\end{algorithm}

\begin{algorithm}
  {\hrulefill}
  \vspace*{-0.4cm}
	\caption{Sequential Monte Carlo Approximate Bayesian Computation Replenishment (RSMC)}
	{\vspace*{-0.25cm}\hrulefill}	
	\label{rsmc}
	\begin{algorithmic}
	  	\STATE Given $N$, $\epsilon_1$, $\epsilon_T$, $c$, $\alpha \in [0,1]$ and $N_\alpha=\lfloor \alpha N \rfloor$,
			\FOR {$i=1$ to $N$} 
		  	\REPEAT 
					\STATE $\mbox{Simulate } \theta_{i} \sim \pi(\theta) \mbox{ and }x\sim f(x|\theta_{i})$
					\STATE $\rho_i =\rho(S(x),S(y))$
				\UNTIL $\rho_i \leq \epsilon_1$
			\ENDFOR
			\STATE Sort $(\theta_i,\rho_i)$ by $\rho_i$
			\STATE Set $\epsilon_{MAX}=\rho_N$
			\WHILE{$\epsilon_{MAX}>\epsilon_T$}
			    \STATE Remove the $N_\alpha$ particles with largest $\rho$
			    \STATE Set $\epsilon_{NEXT}=\rho_{N-N_\alpha}$
			    \STATE Set $i_{acc} =0$
			    \STATE Compute the parameters of the proposal $MCMC$ $q(\cdot,\cdot)$ with the $N-N_\alpha$ particles.
			    \FOR {$j=1$ to $N_\alpha$}
			    	\STATE Simulate $\theta_{N-N_\alpha+j} \sim (\theta_i)_{1\leq i \leq N-N_\alpha}$
			    	\FOR {$k=1$ à $R$} 
			    		\STATE Generate $\theta^* \sim q(\theta^*,\theta_{N-N_\alpha+j})$ et $x^* \sim f(x^*|\theta^*)$
			    		\STATE Generate $u < \mathcal{U}_{[0,1]}$
								\IF{$u \leq 1 \wedge \displaystyle{\frac{\pi(\theta^*)q(\theta_{N-N_\alpha+j},\theta^*)}
								{\pi(\theta_{N-N_\alpha+j})q(\theta^*,\theta_{N-N_\alpha+j})}}\mathds{1}_{\rho(S(x^*),S(y))\leq \epsilon_{NEXT}}$}
									\STATE Set $\theta_{N-N_\alpha+j}=\theta^*$
									\STATE Set $\rho_{N-N_\alpha+j}=\rho(S(x^*),S(y))$ 
									\STATE $i_{acc} \leftarrow i_{acc}+1$
	  						\ENDIF	
				      \ENDFOR 
			    \ENDFOR 
			    \STATE Set $p_{acc}=\displaystyle{\frac{i_{acc}}{RN_\alpha}}$
			    \STATE Set $R=\displaystyle{\frac{\log(c)}{\log(1-p_{acc})}}$
		  \ENDWHILE
	\end{algorithmic}
	{\vspace*{-0.2cm}\hrulefill}		
\end{algorithm}

\begin{algorithm}
  {\hrulefill}
  \vspace*{-0.4cm}
	\caption{Adaptive Sequential Monte Carlo Approximate Bayesian Computation (SMC)}
	{\vspace*{-0.25cm}\hrulefill}	
	\label{smc}
	\begin{algorithmic}
	  \STATE Given $N$, $M$, $\alpha \in [0,1]$, $\epsilon_0 = \infty$, $\epsilon$ and $N_T$,
	  \STATE For $t=0$,
			\FOR {$i=1$ to $N$} 
				\STATE $\mbox{Simulate } \theta_{i}^{(0)} \sim \pi(\theta)$ 
					\FOR {$k=1$ à $M$}
					   \STATE $\mbox{Simulate } X_{(i,k)}^{(0)} \sim f(\cdot|\theta_{i}^{(0)})$
				  \ENDFOR
				\STATE $\mbox{Set } \displaystyle{W_{i}^{(0)}=\frac{1}{N}}$
			\ENDFOR
			\STATE $\mbox{We have } ESS((W_{i}^{(0)}),\epsilon_0)=N \mbox{ where } ESS((W_{i}^{(0)}),\epsilon_0)=\left(\sum_{i=1}^N (W_{i}^{(0)})^2 \right)^{-1}$
			\STATE Set $t=1$
			
		\WHILE {$\epsilon_{t-1}>\epsilon$}
				\STATE Determine $\epsilon_t$ resolving $ESS((W_{i}^{(t)}),\epsilon_t)=\alpha ESS((W_{i}^{(t-1)}),\epsilon_{t-1})$ 
				where $W_{i}^{(t)} \propto W_{i}^{(t-1)} \displaystyle{\frac{\sum_{k=1}^M  \mathds{1}_{A_{\epsilon_{t-1},y}} (X_{(i,k)}^{(t-1)})}{\sum_{k=1}^M  
				\mathds{1}_{A_{\epsilon_{t-1},y}} (X_{(i,k)}^{(t-1)})}}$ et $A_{\epsilon,y}=\left\{x |\,\rho(S(x),S(y))<\epsilon\right\}$
		    \IF {$\epsilon_t < \epsilon$}
		       \STATE $\epsilon_n = \epsilon$ 
		    \ENDIF 	
		    \IF {$ESS((W_{i}^{(t)}),\epsilon_t)<N_T$}
		    	\FOR {$i=1$ to $N$}
		    		\STATE Simulate $(\theta_{(i)}^{(t-1)},X_{(i,1:M)}^{(t-1)})$ in $(\theta_{(j)}^{(t-1)},X_{(j,1:M)}^{(t-1)}) \mbox{ with probabilities } W_{j}^{(t)},\,1 \leq j \leq N$
		    		\STATE Set $W_{i}^{(t)}=\frac{1}{N}$ 
		      \ENDFOR
		    \ENDIF
		    \FOR {$t=1$ to $N$}
		      \IF {$W_{j}^{(t)}>0$} 
						\STATE Generate $\theta^* \sim K(\theta^*|\theta_{(i)}^{(t-1)})$
						\FOR {$k=1$ to $M$}
					  	\STATE $\mbox{Simulate } X_{(*,k)} \sim f(\cdot|\theta^*)$
				    \ENDFOR
						\STATE Generate $u < \mathcal{U}_{[0,1]}$
						\IF{$u \leq 1 \wedge \displaystyle{\frac{\sum_{k=1}^M  \mathds{1}_{A_{\epsilon_{t},y}} (X_{(*,k)})  			
						\pi(\theta^*)K_t(\theta_{(i)}^{(t-1)}|\theta^*)}{\sum_{k=1}^M  \mathds{1}_{A_{\epsilon_{t},y}} (X_{(i,k)}^{(t-1)}) 
						\pi(\theta_{(i)}^{(t-1)})K_t(\theta^*|\theta_{(i)}^{(t-1)})}}$}
							\STATE Set $(\theta_{(i)}^{(t)},X_{(i,1:M)}^{(t)})=(\theta^*,X_{(*,1:M)})$ 
						\ELSE 
							\STATE Set $(\theta_{(i)}^{(t)},X_{(i,1:M)}^{(t)})=(\theta_{(i)}^{(t-1)},X_{(i,1:M)}^{(t-1)})$
						\ENDIF	
					\ENDIF	
			\ENDFOR
		\ENDWHILE	
	\end{algorithmic}
	{\vspace*{-0.2cm}\hrulefill}		
\end{algorithm}

\begin{algorithm}
  {\hrulefill}
  \vspace*{-0.4cm}
	\caption{Adaptive Population Monte Carlo Approximate Bayesian Computation}
	{\vspace*{-0.25cm}\hrulefill}	
	\label{apmc}
	\begin{algorithmic}
	  \STATE Given $N$, $N_\alpha=\lfloor \alpha N \rfloor$ the number of particles to keep at each iteration among the $N$ particles ($\alpha \in [0,1]$) and $p_{acc_{min}}$ the minimal acceptance rate.	
	  \FOR {$t=1$}
			\FOR {$i=1$ to $N$} 
				\STATE Simulate $\theta_{i}^{(0)} \sim \pi(\theta) \mbox{ and } x\sim f(x|\theta_{i}^{(0)})$
				\STATE Set $\rho_{i}^{(0)}=\rho(S(x),S(y))$
        \STATE Set $\displaystyle{w_{i}^{(0)}=1}$	
			\ENDFOR
			\STATE Let $\epsilon_1=Q_{\rho^{(0)}}(\alpha) \mbox{ the first } \alpha \mbox{-quantile of } \rho^{(0)} \mbox{ where } \rho^{(0)}=\left\{\rho_{i}^{(0)}\right\}_{1 \leq i \leq 
			N}$
      \STATE Let $\left\{(\theta_{i}^{(1)},w_{i}^{(1)},\rho_{i}^{(1)})\right\}=\left\{{(\theta_{i}^{(0)},w_{i}^{(0)},\rho_{i}^{(0)})| \rho_{i}^{(0)} \leq \epsilon_1},\ 1 \leq i \leq 
      N\right\}$
      \STATE Take $\sigma_{1}^2$ as twice the weighted empirical variance of $\{(\theta_{i}^{(1)},w_{i}^{(1)})\}_{1 \leq i \leq N_\alpha}$       
      \STATE Set $p_{acc}=1$ 
      	\STATE $t \leftarrow t+1$
    \ENDFOR   
    \WHILE {$p_{acc} > p_{{acc}_{min}}$} 
    	\FOR {$i=N_\alpha+1$ to $N$}
    		\STATE Pick $\theta_i^*$ from $\theta_{j}^{(t-1)} \mbox{ with probability } \frac{w_{j}^{(t-1)}}{\sum_{k=1}^{N_ \alpha} w_{k}^{(t-1)}}$, ${1 \leq j \leq N_\alpha}$
        \STATE Generate $\theta_{i}^{(t-1)}| \theta_i^* \sim \mathcal{N}(\theta_i^*,\sigma_{(t-1)}^2)$ and $x\sim f(x|\theta_{i}^{(t-1)})$
        \STATE Set $\rho_{i}^{(t-1)}=\rho(S(x),S(y))$
        \STATE Set $\displaystyle{w_{i}^{(t-1)}=  \frac{\pi(\theta_{i}^{(t-1)})}{ \sum_{j=1}^{N_\alpha} (w_{j}^{(t-1)}/\sum_{k=1}^{N_ \alpha} w_{k}^{(t-1)}) 
        \sigma_{t-1}^{-1}\varphi(\sigma_{t-1}^{-1} (\theta_{i}^{(t-1)}-\theta_{j}^{(t-1)}))}}$
     	\ENDFOR
     	\STATE Set $p_{acc}=\frac{1}{N-N_\alpha}\sum_{k=N_\alpha+1}^N \mathds{1}_{\rho_{i}^{(t-1)} < \epsilon_{t-1}}$
      \STATE Let $\epsilon_t=Q_{\rho^{(t-1)}}(\alpha) \mbox{ where } \rho^{(t-1)}=\left\{\rho_{i}^{(t-1)}\right\}_{1 \leq i \leq N}$
      \STATE Let $\left\{(\theta_{i}^{(t)},w_{i}^{(t)},\rho_{i}^{(t)})\right\}=\left\{{(\theta_{i}^{(t-1)},w_{i}^{(t-1)},\rho_{i}^{(t-1)})| \rho_{i}^{(t-1)} \leq \epsilon_t},\ 1 
      \leq i \leq N\right\}$ 
      \STATE Take $\sigma_{t}^2$ as twice the weighted empirical variance of $\{(\theta_{i}^{(t)},w_{i}^{(t)})\}_{1 \leq i \leq N_\alpha}$
     	\STATE $t \leftarrow t+1$
		\ENDWHILE	
		\STATE Where $\forall u \in [0,1]$ and $X=\{x_1,...,x_n\},\,Q_X(u)=\inf\{x \in X | F_X(x)\geq u\}$ and $F_X(x)=\frac{1}{n} \sum_{k=1}^n \mathds{1}_{x_k \leq x}$.
    \STATE Where $\varphi(x)=\frac{1}{\sqrt{2\pi}} e^{-\frac{x^2}{2}}$	
 	\end{algorithmic}
	{\vspace*{-0.2cm}\hrulefill}		
\end{algorithm}

\clearpage
\newpage
\section*{Appendix B: Proof that the algorithm stops}

We know that there exists $\epsilon_\infty > 0$ such that $\epsilon_t \underset{t\to+\infty}{\longrightarrow} \epsilon_\infty$ because, by construction of the algorithm  $(\epsilon_t)$ is a positive decreasing sequence and it is bounded by 0.   

\noindent For each $\theta \in \Theta$, we consider the distance $(\rho(x,y)| \theta)$ as a random variable $\rho(\theta)$. Let $f_{\rho(\theta)}$ be the probability density function of $\rho(\theta)$.

\noindent The probability $\mathbb{P}[\rho(\theta) \geq \epsilon_{t}] $ that the drawn distance associated to parameter $\theta$ is higher than the current tolerance $\epsilon_{t}$ satisfies:

\noindent 

	$$\begin{array}{ll}
\mathbb{P}[\rho(\theta) \geq \epsilon_{t}] & = 1 -  \mathbb{P}[(\rho(\theta) < \epsilon_{t}]\\
																					 & = 1 - \int_{\epsilon_\infty}^{\epsilon_{t}} f_{\rho(\theta)}(x) d_x\\

  \end{array}$$

\noindent  We define: 

$$\mathbb{P}_{max}=\sup_{\theta \in \Theta}\left\{\sup_{x \in \mathbb{R}^+}\left\{f_{\rho(\theta)}(x)\right\}\right\}$$

\noindent We have:

	$$
\mathbb{P}[\rho(\theta) \geq \epsilon_{t}] \geq 1 -\mathbb{P}_{max}(\epsilon_{t}-\epsilon_\infty)
$$

\noindent The $N-N_\alpha$ particles are independent and identically distributed from $\pi_{t+1}$ the density defined by the algorithm, hence the probability $ \mathbb{P}[p_{acc}(t+1)=0] $ that no particle is accepted at step $t+1$ is such that:

\noindent 

			$$
            \mathbb{P}[p_{acc}(t+1)=0] \geq  \left( 1 -  \mathbb{P}_{max}(\epsilon_{t}-\epsilon_\infty)\right)^{N-N_\alpha}
 $$

\noindent If $\mathbb{P}_{max} < +\infty$, because $\epsilon_t - \epsilon_{\infty}\underset{t\to+\infty}{\longrightarrow} 0$, we have:

 $$\mathbb{P}[p_{acc}(t+1)=0]\underset{t\to+\infty}{\longrightarrow} 1$$

\noindent We can conclude that $p_{acc}(t)$ converges in probability towards 0 if $\mathbb{P}_{max} < +\infty$. This ensures that the algorithm stops, whatever the chosen value of $p_{acc_{min}}$.

\end{document}